\documentclass[reqno]{amsart}

\usepackage{amsthm, amssymb, latexsym, tabls, hhline,lscape}
\usepackage{color}

\usepackage[margin=1.5in]{geometry}

\newtheorem{theorem}{Theorem}[section]
\newtheorem{proposition}[theorem]{Proposition}

\newtheorem{cor}[theorem]{Corollary}
\newtheorem{conjecture}[theorem]{Conjecture}
\newtheorem{lemma}[theorem]{Lemma}
\theoremstyle{definition}

\newtheorem{example}[theorem]{Example}

\theoremstyle{remark}
\newtheorem{remark}[theorem]{Remark}
\numberwithin{equation}{section}
\newcommand{\Sym}{\ensuremath{\mathit{Sym}}}
\newcommand{\QSym}{\ensuremath{\mathit{QSym}}}
\newcommand{\NSym}{\ensuremath{\mathit{NSym}}}
\newcommand{\ra}{\rightarrow}
\newcommand{\ot}{\ensuremath{\otimes}}

\newcommand{\scal}[1]{\langle #1 \rangle}

\newcommand{\Z}{\mathbb Z}
\newcommand{\N}{\mathbb N}
\newcommand{\C}{\mathbb C}

\newcommand{\s}{\mathfrak{S}}
\newcommand{\Des}{\mathrm{Des}}

\makeindex

\title[Combinatorial Hopf algebras and Towers of Algebras]
{\bf Combinatorial Hopf algebras and Towers of Algebras -- \\
Dimension, Quantization and Functorality}

\author{N. Bergeron}\address[Nantel Bergeron]
{Department of Mathematics and Statistics\\ York  University\\ To\-ron\-to, Ontario M3J 1P3\\ CANADA}
\email{bergeron@mathstat.yorku.ca}
\urladdr{http://www.math.yorku.ca/bergeron}

 \author{T. Lam}\address[Thomas Lam]
{Department of Mathematics \\ University of Michigan \\
Ann Arbor, MI 48109
\\ USA}
 \email{tfylam@umich.edu}
 \urladdr{http://www.math.lsa.umich.edu/\~{ }tfylam}

 \author{H. Li}\address[Huilan Li]
 {Department of Mathematics\\ Drexel University\\
 Philadelphia\\ PA 19104 \\ USA}
 \email{huilan.li@gmail.com}
  \urladdr{http://www.math.drexel.edu/\~{}huilan}

\date{\today}

\thanks{N.B. is supported in part by CRC and NSERC}
\thanks{T.L. is partially supported by NSF grants DMS-0600677 and DMS-0652641}
\thanks{H.L. is supported in part by CRC, NSERC and NSF 06-580}



\begin{document}
\maketitle


\begin{abstract}
Bergeron and Li have introduced a set of axioms which guarantee that
the Grothendieck groups of a tower of algebras
$\bigoplus_{n\ge0}A_n$ can be a pair of graded dual Hopf algebras.
Hivert and Nzeutchap, and independently Lam and Shimozono
constructed dual graded graphs from primitive elements in Hopf
algebras.  In this paper we apply the composition of these
constructions to towers of algebras. We show that if a tower
$\bigoplus_{n\ge0}A_n$ gives rise to graded dual Hopf algebras, then
$\dim(A_n)=r^nn!$ where $r = \dim(A_1)$. In the case $r=1$ we give a
conjectural classification. We then investigate a quantum version of
the main theorem. We conclude with some open problems and a
categorification of these constructions.
\end{abstract}

\section{Introduction}
\setcounter{equation}{0}

This paper\footnote{A summary of an earlier shorter version of this
paper appeared in \cite{BLL}.} is concerned with the interplay
between towers of associative algebras, pairs of graded dual
combinatorial Hopf algebras, and dual graded graphs. Our point of
departure is the study of the composition of two constructions: (1)
the construction of dual Hopf algebras from towers of algebras
satisfying some axioms, due to Bergeron and Li \cite{BL}; and (2)
the construction of dual graded graphs from primitive elements in
dual Hopf algebras, discovered independently by Hivert and Nzeutchap
\cite{HN}, and Lam and Shimozono:
\begin{equation}\label{eq:comp}
\text{tower of algebras} \longrightarrow \text{combinatorial Hopf
algebra} \longrightarrow \text{dual graded graph}
\end{equation}

\subsection{From towers of algebras to combinatorial Hopf algebras}
A {\it tower} $\bigoplus_{n \geq 0} A_n$ of associative algebras is a
sequence of associative algebras $A_n$ over $\C$, equipped with an
external multiplication $\rho_{m,n}: A_m \otimes A_n \to A_{m+n}$
satisfying a number of axioms (see Section \ref{sec:BL}).  A {\it
combinatorial Hopf algebra} is a graded, connected Hopf algebra with
a distinguished basis such that all product and coproduct structure
constants are nonnegative.  The Grothendieck groups $G(A)$ and
$K(A)$ of certain towers of algebras give rise to a pair of graded
dual combinatorial Hopf algebras: the product and coproduct
structures come from induction and restriction of modules, while the
distinguished bases comes from the classes of simple modules, and of
indecomposable projective modules.  Our notion of combinatorial Hopf
algebra is related to, but different from the one in \cite{ABS} (see
Section \ref{sec:further}).

\subsection{From combinatorial Hopf algebras to dual graded graphs}
The notion of a pair $(\Gamma,\Gamma')$ of {\it dual graded graphs}
(see Section \ref{sec:HNLS}) was introduced by Fomin \cite{F} (see
also \cite{S}) to encode the enumerative properties of the
Robinson-Schensted correspondence and its generalizations.  These
are pairs of graded graphs equipped with linear operators $D, U$
satisfying the relation $DU - UD = rI$, where $r$ is a nonnegative
integer, and $I$ is the identity.  The second arrow of
(\ref{eq:comp}) is obtained by using (some of the) structure
constants of a combinatorial Hopf algebra as edge multiplicities for
a graph.  This construction depends on some choices, but for
combinatorial Hopf algebras arising from towers of algebras there is
an extremely natural choice (see (\ref{eq:ab})).

\subsection{The dimension theorem}
By relating the dimensions of simple modules and indecomposable
projective modules to paths in dual graded graphs via
(\ref{eq:comp}), we obtain the following result, showing that towers
of algebras giving rise to combinatorial Hopf algebras are much more
rigid than they first appear.

\begin{theorem} \label{dimthm}
Let $A=\bigoplus_{n\geq0}A_{n}$ be a tower of algebras. If the
Grothendieck groups of $A$ form a pair of graded dual Hopf algebras,
then $\dim(A_n)=r^nn!$ where $r=\dim(A_1)$.
\end{theorem}

The notion of ``forming a pair of graded dual Hopf algebras'' is
made precise in Section \ref{sec:BL}.  The number $r^n n!$ counts
certain paths in a pair of dual graded graphs. It is also the number
of permutation matrices with entries in a finite group of size $r$.

\subsection{Symmetric groups, symmetric functions, and Young's graph}
The fundamental example of all three classes of objects, arises from
the representation theory of symmetric groups and the theory of
symmetric functions.  The symmetric group algebras give rise to a
tower $\bigoplus_{n\geq0}\mathbb{C}\mathfrak{S}_{n}$ of algebras,
the Grothendieck groups of which form a graded self-dual
combinatorial Hopf algebra (see \cite{Bur} and \cite{AVZ}), which
can be identified with the ring \Sym\ of symmetric functions. The
Hopf structure of \Sym\ was studied by Geissinger \cite{LG}.  The
corresponding dual graded graph is the (self-dual) {\it Young graph}
of partitions, which is the motivating example \cite{F,S} of dual
graded graphs. Indeed Young's graph can be obtained from Young's
branching rule for the symmetric group, or equivalently the Pieri
rule for symmetric functions.

\subsection{Towards a classification?}
In recent years it has been shown that other graded dual Hopf
algebras can be obtained from towers of algebras.  In~\cite{MR}
Malvenuto and Reutenauer establish the duality between the Hopf
algebra \NSym\ of noncommutative symmetric functions and the Hopf
algebra \QSym\ of quasi-symmetric functions. Krob and
Thibon~\cite{KT} then showed that this duality can be interpreted as
the duality of the Grothendieck groups associated with
$\bigoplus_{n\geq0}H_{n}(0)$ the tower of Hecke algebras at $q = 0$.
For more examples, see~\cite{BHT,HNT,ANS}.

It is very tempting, as suggested by J. Y. Thibon, to classify all
combinatorial Hopf algebras which arise as Grothendieck groups
associated with a tower of algebras $\bigoplus_{n\geq0}A_{n}$. The
list of axioms given by the first and last author in \cite{BL}
guarantees that the Grothendieck groups of a tower of algebras form
a pair of graded dual Hopf algebras.  This list of axioms is far
from a classification.

The rigidity proved in Theorem \ref{dimthm} suggests however that
there may be a classification theorem for towers of algebras which
give rise to combinatorial Hopf algebras. For the case $r=1$, we
give a conjectural classification in Section~\ref{sec:HC}, which
includes symmetric group algebras, 0-Hecke algebras, nilCoxeter
algebras (studied by Khovanov \cite{Kho}), and the infinite family
of Hecke algebras at a root of unity (see \cite{LLT}).  The
representation theory of Hecke algebras at roots of unity is a very
active area of research.

In general, Theorem \ref{dimthm} suggests that to perform the
inverse constructions of the arrows in (\ref{eq:comp}) one should
study algebras related to symmetric groups (or wreath products of
symmetric groups).  There are many combinatorial Hopf algebras for
which one may attempt to perform the inverse construction, but there
are even more dual graded graphs.  The general construction of
\cite{LamSKac} produces dual graded graphs from Bruhat orders of
Weyl groups of Kac-Moody algebras and it is unclear whether there
are Hopf algebras, or towers of algebras giving rise to these
graphs.

\subsection{Quantization and categorification}
Our work should also be put into the context of the general notion
of categorification: Theorem \ref{dimthm} provides a condition on
when a Hopf algebra can be categorified by the Grothendieck groups
of a tower of algebras.  The idea of categorifying Hopf algebras, in
particular quantum groups, has been around for some time.  For
example, Crane and Frenkel \cite{CF} introduced a notion of a {\it
Hopf-category} in the context of four-dimensional topological
quantum field theory. We remark that quantum groups are not graded,
and so do not fit into the class of Hopf algebras that we consider.

In the last two sections, we turn our attention to some
generalizations.  In Section~\ref{Sec:Q}, we give a quantum version
of Theorem~\ref{dimthm}: replacing towers of algebras with filtered
towers of algebras, Hopf algebras with $q$-twisted Hopf algebras
\cite{LZ}, and dual graded graphs with quantized dual graded graphs
\cite{Lam}.

In Section \ref{sec:further}, we relate our work to the
combinatorial Hopf algebras of Aguiar, Bergeron and
Sottile~\cite{ABS}.  We also discuss different notions of towers of
algebras, and describe how to categorify the constructions in
(\ref{eq:comp}), in particular making the arrows into functors.

{\bf Acknowledgements.}  We thank Susumu Ariki for comments.  Lam would like to thank Mark Shimozono for the
collaboration which led to the line of thinking in this paper.


\section{From towers of algebras to combinatorial Hopf algebras}
\label{sec:BL} We recall here the work of Bergeron and Li~\cite{BL}
on towers of algebras.  
For $B$ an arbitrary algebra we denote by $
_{B}\mbox{mod}$, the category of all finitely generated left
$B$-modules, and by $\mathcal{P}(B)$, the category of all finitely
generated projective left $B$-modules. For some category
$\mathcal{C}$ of left $B$-modules ($_{B}\mbox{mod}$ or
$\mathcal{P}(B)$) let $\mathbf{F}$ be the free abelian group
generated by the symbols $(M)$, one for each isomorphism class of
modules $M$ in $\mathcal{C}$. Let $\mathbf{F_{0}}$ be the subgroup
of $\mathbf{F}$ generated by all expressions $(M)-(L)-(N)$ one for
each exact sequence
              $$0\rightarrow L\rightarrow M\rightarrow N\rightarrow0$$
in $\mathcal{C}$. The
\textit{Grothendieck group} ${\mathcal K}_{0}(\mathcal{C})$ of the category
$\mathcal{C}$ is defined by the quotient $\mathbf{F}/\mathbf{F_{0}}$, an abelian additive group.
For $M\in\mathcal{C}$, we denote by $[M]$ its image in ${\mathcal K}_{0}(\mathcal{C})$.
We then set
             $$G_{0}(B)={\mathcal K}_{0}(_{B}\mbox{mod}) \quad\hbox{and}\quad
            K_{0}(B)={\mathcal K}_{0}(\mathcal{P}(B)).$$

For $B$ a finite-dimensional algebra over a field $K$, let
$\{V_{1},\cdots, V_{s}\}$ be a complete list of nonisomorphic
simple $B$-modules. The projective covers $\{P_{1},\cdots,
P_{s}\}$ of the simple modules $V_i$'s  is a complete list of nonisomorphic indecomposable
projective $B$-modules. Then  $G_{0}(B)=\bigoplus^s_{i=1}\mathbb{Z}[V_{i}]$ and
$ K_{0}(B)=\bigoplus^s_{i=1}\mathbb{Z}[P_{i}]$.

Let $\varphi\colon B\to A$ be an injection of algebras preserving
unities, and let $M$ be a (left) $A$-module and $N$ a (left)
$B$-module. The \textit{induction} of $N$ from $B$ to $A$ is
Ind$^{A}_{B}N=A\otimes_{\varphi}N$, the (left) $A$-module $A\otimes
N$ modulo the relations $a\otimes bn\equiv a\varphi(b)\otimes n$,
and the \textit{restriction} of $M$ from $A$ to $B$ is
Res$^{A}_{B}M=\mbox{Hom}_{A}(A,M)$, the (left) $B$-module with the
$B$-action
 defined by $bf(a)=f(a\varphi(b))$.

Let $A=\bigoplus_{n\geq 0}A_{n}$ be a graded algebra over
$\mathbb{C}$ with multiplication $\rho\colon A\otimes A\to A$.
 Bergeron and Li studied five axioms for $A$ (we refer to \cite{BL} for full details):

\smallskip
\noindent(1) For each $n\ge 0$, $A_{n}$ is a finite-dimensional
algebra by itself with (internal) multiplication $\mu_n\colon
A_n\otimes A_n\to A_n$  and unit $1_n$. $A_{0}\cong \mathbb{C}$.

\smallskip
\noindent(2) The (external) multiplication $\rho_{m,n}:
A_{m}\otimes A_{n}\rightarrow A_{m+n}$ is an injective
homomorphism of algebras, for all $m\mbox{ and }n$ (sending
$1_{m}\otimes 1_{n}$ to $1_{m+n}$ ).

\smallskip
\noindent(3) $A_{m+n}$ is a two-sided projective $A_{m}\otimes
A_{n}$-module with the action defined by $a\cdot(b\ot
c)=a\rho_{m,n}(b\ot c)\mbox{ and }(b\ot c)\cdot a=\rho_{m,n}(b\ot
c)a$, for all $m,n\geq0,\ a\in A_{m+n},\ b\in A_{m},\ c\in
A_{n}\mbox{ and }m,n\geq0$.

\smallskip
\noindent(4) A relation between the decomposition of $A_{n+m}$ as a
left $A_m\ot A_n$-module and as a right $A_m\ot A_n$-module holds.

\smallskip
\noindent(5) An analogue of Mackey's formula relating induction and restriction of modules holds.

\medskip
We say here that $A = \bigoplus_{n\geq 0}A_{n}$ is a {\it tower of
algebras} if it satisfies Conditions (1), (2) and (3).

Condition (1) guarantees that the Grothendieck groups
$G(A)=\bigoplus_{n\ge 0} G_0(A_n)$ and $K(A)=\bigoplus_{n\ge 0}
K_0(A_n)$ are graded connected.  Conditions (2) and (3) ensure that
induction and restriction are well defined on $G(A)$ and $K(A)$,
defining a multiplication and comultiplication, as follows.  For
$[M]\in G_0(A_m)$ (or $K_0(A_m)$) and $[N]\in G_0(A_n)$ (or
$K_0(A_n)$) we let
 $$[M] [N]=\left[\mbox{Ind}^{A_{m+n}}_{A_{m}\otimes
A_{n}}M\otimes N\right]  \qquad\hbox{and}\qquad
\Delta([N])=\sum_{k+l=n} \left[\mbox{Res}^{A_{k+l}}_{A_{k}\otimes
A_{l}} N\right].$$ The pairing between $K(A)$ and $G(A)$ is given by
$\scal{\,\,,\,}: K(A)\times G(A)\ra \mathbb{\mathbb{Z}} $
where
        $$\scal{[P],[M]}=\left\{\begin{array}{ll}
                       \mbox{dim}_{K}\big(\mbox{Hom}_{A_{n}}(P,M)\big)
                       & \mbox{if }[P]\in K_{0}(A_{n})
                       \mbox{ and }[M]\in G_{0}(A_{n}),\\
                       0 & \mbox{otherwise.}
                       \end{array}
                \right.$$

Thus with (only) Conditions (1), (2), and (3), $G(A)$ and $K(A)$ are
dual free $\Z$-modules both endowed with a multiplication and
comultiplication.  
\begin{theorem}[Bergeron and Li \cite{BL}]
If a graded algebra $A=\bigoplus_{n\geq 0}A_{n}$ over $\mathbb{C}$
satisfies Conditions (1)-(5), then $G(A)$ and $K(A)$ are graded dual
Hopf algebras.
\end{theorem}
In particular Theorem \ref{dimthm} applies to graded algebras which
satisfy Conditions (1)-(5).  Note that the dual Hopf algebras $G(A)$
and $K(A)$ come with distinguished bases consisting of the
isomorphism classes of simple and indecomposable projective modules.

\medskip
\section{From combinatorial Hopf algebras to dual graded graphs}
\label{sec:HNLS} This section recounts work of Fomin \cite{F},
  Hivert and Nzeutchap \cite{HN}, and Lam and Shimozono.  A {\it
graded graph} $\Gamma = (V,E,h,m)$ consists of a set of vertices
$V$, a set of (directed) edges $E \subset V \times V$, a height
function $h: V \to \{0,1,\ldots\}$ and an edge multiplicity function
$m: V \times V \to \{0,1,\ldots\}$. If $(v,u) \in E$ is an edge, then
we must have $h(u) = h(v) + 1$. The multiplicity function determines
the edge set: $(v,u) \in E$ if and only if $m(v,u) \neq 0$.  We
assume always that there is a single vertex $v_0$ of height 0.

Let $\Z V = \bigoplus_{v \in V} \Z \cdot v$ be the free $\Z$-module
generated by the vertex set.  Given a graded graph $\Gamma = (V, E,
h, m)$ we define up and down operators $U,D \colon\Z V \to \Z V$ by
$$
U_\Gamma(v) = \sum_{u \in V}m(v,u)\, u \ \ \ D_{\Gamma}(v) = \sum_{u
\in V}m(u,v)\, u
$$
and extending by linearity over $\Z$.  We will assume that $\Gamma$
is locally-finite, so that these operators are well defined.  A pair
$(\Gamma, \Gamma')$ of graded graphs with the same vertex set $V$
and height function $h$ is called {\it dual} with {\it differential
coefficient} $r$ if 
$$ D_{\Gamma'} U_{\Gamma} - U_\Gamma D_{\Gamma'} =
r\,{\rm Id}.
$$
We shall need the following result of Fomin.  For a graded graph
$\Gamma$, let $f^v_{\Gamma}$ denote the number of paths from $v_0$
to $v$, where for two vertices $w, u \in V$, we think that there are
$m(w,u)$ edges connecting $w$ to $u$.
\begin{theorem}[Fomin~\cite{F}] \label{thm:F} Let $(\Gamma,\Gamma')$ be a pair of dual graded graphs
with differential coefficient $r$. Then
$$r^n n!=\sum_{v \colon h(v)=n}f^v_{\Gamma}f^v_{\Gamma'}.$$

\end{theorem}

Let $H_\bullet = \bigoplus_{n \geq 0} H_n$ and $H^\bullet =
\bigoplus_{n \geq 0} H^n$ be a pair of graded dual Hopf algebras over $\Z$
with respect to the pairing $\scal{\,.\,,.\,}: H_\bullet \times H^\bullet
\to \Z$. We assume that we are given dual sets of homogeneous free
$\Z$-module generators $\{p_\lambda \in H_\bullet\}_{\lambda \in
\Lambda}$ and $\{s_\lambda \in H^\bullet\}_{\lambda \in \Lambda}$,
such that all structure constants are non-negative integers.  We
also assume that $\dim(H_i) = \dim(H^i) < \infty$ for each $i \geq
0$ and $\dim(H_0) = \dim(H^0) = 1$, so that $H_0$ and $H^0$ are
spanned by distinguished elements the unit $1$.  Let us suppose we
are given non-zero homogeneous elements $\alpha \in H_1$ and $\beta
\in H^1$ of degree 1 such that $\alpha p_{\mu}$ (resp. $\beta s_{\mu}$)
is a linear combination of $\{p_{\lambda}\}$ (resp. $\{s_{\lambda}\}$) with nonnegative integer
coefficients for any $\mu\in\Lambda$.

We now define a graded graph $\Gamma(\beta) = (V, E, h, m)$ where
$V=\{s_\lambda\}_{\lambda \in \Lambda}$ and $h\colon V\to{\mathbb
Z}$ is defined by $h(s_\lambda)=\deg(s_\lambda)$.  The map $m\colon
V\times V\to{\mathbb Z}$ is defined by
$$
m(s_\lambda,s_\mu)= \scal{p_\mu,\beta s_\lambda} =
\scal{\Delta(p_\mu), \beta\otimes s_\lambda}
$$
and $E$ is determined by $m$.  The grading of $\Gamma(\beta)$
follows from the assumption that $\beta$ has degree 1.  Similarly,
we define a graded graph $\Gamma'(\alpha)=(V',E',h',m')$ where $V' =
V$, $h' = h$, and
$$
m'(s_\lambda,s_\mu)= \scal{\alpha \, p_\lambda, s_\mu} = \scal{
\alpha \ot  p_\lambda, \Delta(s_\mu)}.$$  The following theorem is
due independently to Hivert and Nzeutchap \cite{HN} and Lam and
Shimozono (unpublished).

\begin{theorem}\label{thm:HNLS}
The graded graphs $\Gamma = \Gamma(\beta)$ and $\Gamma' =
\Gamma'(\alpha)$ form a pair of dual graded graphs with {\it
differential coefficient} $\scal{\alpha,\beta}$.
\end{theorem}
\begin{proof}
We identify $\Z V$ with $H^\bullet$ and note that $U_\Gamma(x) =
\beta \, x$ where $x \in H^\bullet$ and we use the multiplication in
$H^\bullet$. Also,
\begin{align*}
D_{\Gamma'}(x) &= \sum_{\mu \in \Lambda} \scal{\alpha \ot p_\mu,
\Delta x} \, s_\mu = \sum \scal{\alpha,x^{(1)}}\, x^{(2)}.
\end{align*}
where $\Delta x = \sum x^{(1)} \otimes x^{(2)}$. Now observe that by
our hypotheses on the degree of $\alpha$ and $\beta$ they are
primitive elements: $\Delta\alpha = 1 \otimes \alpha + \alpha
\otimes 1$ and $\Delta\beta = 1 \otimes \beta + \beta \otimes 1$. We
first calculate
\begin{align*}
\scal{\alpha,\beta \, x} =\scal{\Delta\alpha,\beta \otimes x}
=\scal{1,\beta}\scal{\alpha,x}+\scal{\alpha,\beta}\scal{1,x}
=\scal{\alpha,\beta}\scal{1,x}
\end{align*}
and then compute
\begin{align*}
D_{\Gamma'} U_{\Gamma}(x) & = D_{\Gamma'}(\beta x)\\
& = \sum \left(\scal{\alpha,\beta \, x^{(1)}}\,x^{(2)} +
\scal{\alpha,x^{(1)}}\,\beta\,x^{(2)}\right) \\
&= \scal{\alpha,\beta}x +U_{\Gamma}D_{\Gamma'}(x)
\end{align*}
where to obtain $\scal{\alpha,\beta}x$ in the last line we use
$\Delta x = 1 \ot x + \text{terms of other degrees}.$
\end{proof}

\section{Proof of Theorem \ref{dimthm}}\label{sec:proof}
We are given a graded algebra $A=\bigoplus_{n\geq 0}A_{n}$ over
$\mathbb{C}$ with multiplication $\rho$ satisfying Conditions (1),
(2) and (3).  Moreover we assume that the two Grothendieck groups
$G(A)$ and $K(A)$ form a pair of graded dual Hopf algebras as in
Section \ref{sec:BL}.  Under these assumptions we show that
$$ \dim(A_n)=r^n n! $$ where $r=\dim(A_1)$.

Let $H^{\bullet}=G(A)$ and $H_{\bullet}=K(A)$.  Let $\{s^{(1)}_1 =
[S^{(1)}_1],\ldots,s^{(1)}_t = [S^{(1)}_t]\}$ and $\{p^{(1)}_1=
[P^{(1)}_1],\ldots,p^{(1)}_t = [P^{(1)}_t]\}$ denote the isomorphism
classes of simple and indecomposable projective $A_1$-modules, so
that $H^1=\bigoplus_{i=1}^{t}\mathbb{Z}s^{(1)}_i$ and
$H_1=\bigoplus_{i=1}^{t}\mathbb{Z}p^{(1)}_i$.  Define
$a_i=\dim(S^{(1)}_i)$ and $b_i=\dim(P^{(1)}_i)$ for $1\leq i\leq t$.
We set for the remainder of this paper
  \begin{equation}\label{eq:ab}
  \alpha=\sum_{i=1}^ta_ip^{(1)}_i\in H_{1}\qquad \hbox{and}\qquad \beta=\sum_{i=1}^tb_is^{(1)}_i\in H^{1}.
  \end{equation}

Since $A_0\cong \mathbb{C}$, we let $s^{(0)}_1$ (respectively,
$p^{(0)}_1$) be the unique simple (respectively, indecomposable
projective) module representative in $H^0$ (respectively, $ H_0$).
Similarly, let $\{s^{(n)}_i = [S^{(n)}_i]\}$ be all isomorphism
classes of simple $A_n$-modules and $\{p^{(n)}_i = [P^{(n)}_i]\}$ be
all isomorphism classes of indecomposable projective $A_n$-modules.
The sets $\bigcup_{n\ge 0}\{s^{(n)}_i\}$ and $\bigcup_{n\ge
0}\{p^{(n)}_i\}$ form dual free $\Z$-module bases of $H^{\bullet}$
and $H_{\bullet}$.

Now define $\Gamma = \Gamma(\beta)$ and $\Gamma' = \Gamma'(\alpha)$
as in Section \ref{sec:HNLS}.

\begin{lemma}\label{lem:dimpath}
The numbers of paths from $s^{(0)}_1$ to ${s^{(n)}_j}$ in $\Gamma$ and $\Gamma'$ are
$$
f^{s^{(n)}_j}_{\Gamma} = \dim P^{(n)}_j \ \ \text{and} \ \
f^{s^{(n)}_j}_{\Gamma'} = \dim S^{(n)}_j.
$$
\end{lemma}
\begin{proof}
From the definition of $\beta$ in (\ref{eq:ab}),
$$
m(s^{(n-1)}_i,s^{(n)}_j)=\sum_{l=1}^t b_lc_l,
$$
where $c_l$ is the number of copies of the indecomposable projective
module $P^{(1)}_l\otimes P^{(n-1)}_i$ as a summand in
Res$_{A_{1}\otimes A_{n-1}}^{A_n}P^{(n)}_j$.  Note that $s^{(0)}_1$
is the unit of $H^{\bullet}$ and $m(s^{(0)}_1,s^{(1)}_i)=b_i=\dim
P^{(1)}_i$ for all $1\leq i\leq t$. The dimension of an
indecomposable projective module $P^{(n)}_j$ is given by
  $$\dim P^{(n)}_j = \sum_{i,l} c_l \dim  \left(P^{(1)}_l\otimes P^{(n-1)}_i\right)= \sum_i m(s^{(n-1)}_i,s^{(n)}_j)\dim P^{(n-1)}_i.
  $$
By induction on $n$,  we deduce that $\dim P^{(n)}_j$ is the number
of paths from $s^{(0)}_1$ to $s^{(n)}_j$ in $\Gamma$.  The claim for
$\Gamma'$ is similar.
\end{proof}

For any  finite dimensional algebra $B$ let
$\{S_{\lambda}\}_{\lambda}$ be a complete set of simple $B$-modules.
For each $\lambda$ let $P_\lambda$ be the projective cover of
$S_\lambda$. It is well known (see~\cite{CR}) that we can find
minimal idempotents $\{e_i\}$ such that $B=\bigoplus Be_i$ where
each $Be_i$ is isomorphic to a $P_\lambda$. Moreover, the quotient
of $B$ by its radical shows that the multiplicity of $P_\lambda$ in
$B$ is equal to $\dim S_\lambda$. This implies the following lemma.

\begin{lemma}\label{dimreplem} Let $B$ be a finite dimensional  algebra and $\{S_{\lambda}\}_{\lambda}$ be a complete set of simple $B$-modules.
$$\dim B=\sum_{\lambda} (\dim P_{\lambda})(\dim S_{\lambda}),$$
where $P_{\lambda}$ is the  projective cover of $S_{\lambda}$.
\end{lemma}

By Lemma \ref{dimreplem}, $r= \sum_{i=1}^t a_ib_i=
\scal{\alpha,\beta}$.  By Theorem \ref{thm:HNLS} we apply
Theorem \ref{thm:F} to $(\Gamma,\Gamma')$.  Using Lemma
\ref{dimreplem} and Lemma \ref{lem:dimpath}, Theorem \ref{thm:F}
says
$$
\dim(A_n) = \sum_i (\dim P^{(n)}_{i})(\dim S^{(n)}_{i}) = \sum_i
f_\Gamma^{s_i^{(n)}}\,f_{\Gamma'}^{s_i^{(n)}} = r^n n!.
$$

\begin{remark} If the tower consists of semisimple algebras $A_i$,
then $\Gamma = \Gamma'$. So we obtain a {\it self-dual graded graph} $\Gamma$.  In
this case the graph would be a weighted version of a {\it
differential poset} in the sense of Stanley \cite{S}.  If
furthermore the branching of irreducible modules from $A_n$ to
$A_{1} \otimes A_{n-1}$ is multiplicity free, then we get a true
differential poset.
\end{remark}

\begin{remark} The Hopf algebras $H^\bullet$ and $H_\bullet$ are not
in general either commutative or co-commutative.  Thus in the
definitions of Section \ref{sec:HNLS} we could have obtained a
different pair of dual graded graphs by setting $m(s_\lambda,s_\mu)=
\scal{p_\mu, s_\lambda\, \beta}$  or $m'(s_\lambda,s_\mu)=
\scal{p_\lambda \, \alpha, s_\mu}$.
\end{remark}

\section{Examples 
}\label{sec:ex} In this section we explain four examples of the
constructions in Sections \ref{sec:BL}-\ref{sec:proof}, all with $r
= \dim(A_1)=1$.

\def\J{{\mathcal J}}

\begin{center}
{\small
\begin{tabular}{|c|c|c|c|c|}
\hline
Tower of algebras $A$ & $K(A)$ & $G(A)$ & $\Gamma$ & $\Gamma'$ \\
\hline
\hline
$\s_n$-group algebras & $\Sym$ & $\Sym$ & Young's graph & Young's graph \\
\hline
NilCoxeter algebras & $\Z[x]$ & $\Z[x, x^2/2, x^3/3!,\ldots]$ & weighted chain & chain \\
\hline
0-Hecke algebras & $\NSym$ & $\QSym$ & BinWord graph & lifted binary tree \\
\hline
Hecke algebras at $\sqrt[r]{1}$& $(\J^{(r)})^\perp$ & $\Sym/\J^{(r)}$ & ??? & ??? \\
\hline
\end{tabular}
}
\end{center}
\subsection{Symmetric group algebras}
Let $A = \bigoplus_{n \geq 0} \C\s_n$ be the tower of symmetric group algebras.  Since $\C\s_n$ is semisimple, $K(A) = G(A)$.  Indeed both $K(A)$ and $G(A)$ can be identified with the (self-dual) Hopf algebra $\Sym$ of symmetric functions, and the classes of the simple modules and the indecomposable projective modules are identified with the Schur functions $s_\lambda$.  The corresponding self-dual graded graph is Young's lattice of partitions.  We refer the reader to \cite{AVZ} for further details of this well-known example.

\subsection{NilCoxeter algebras}\label{ssec:nilCox}
The {\it nilCoxeter algebra} $N_n$ is the unital algebra over $\C$ generated by $T_1,T_2,\ldots, T_{n-1}$ with relations
\begin{align*}
T_i^2 &= 0 \\
T_iT_j &= T_jT_i  & \mbox{for $|i-j| > 1$} \\
T_i T_{i+1} T_i &= T_{i+1} T_i T_{i+1}.
\end{align*}
It has a basis $\{T_w \mid w \in \s_n\}$ labeled by permutations of $\{1,2,\ldots,n\}$, where $T_w = T_{i_1}T_{i_2}\cdots T_{i_\ell}$ if $w = s_{i_1}s_{i_2} \cdots s_{i_\ell}$ is a reduced factorization of $w$.  An explicit realization of this algebra is obtained by {\it divided difference operators}.  The external multiplication $N_i \otimes N_j \to N_{i+j}$ is defined in the same way as for symmetric group algebras.  The representation theory of the tower $N= \bigoplus_{n \geq 0} N_n$ was worked out by Khovanov \cite{Kho}.

The unique simple module (up to isomorphism) $S_n$ of $N_n$ has dimension 1, with projective cover $P_n = N_n$.  Then $K(N) = \Z[x]$ with $[P_n]  = x^n$ and $G(N) = \Z[x,x^2/2, x^3/3!, \ldots]$ with $[S_n] = x^n/n!$. Here the algebra $\Z[x,x^2/2, x^3/3!, \ldots]$ is the free divided powers algebra on one generator over $\Z$. The coproduct is given by $\Delta(x) = 1 \otimes x + x \otimes 1$ for both $K(N)$ and $G(N)$.  The graph $\Gamma'$ is a chain, with vertices $\{0,1,2,\ldots\}$ and multiplicities $m(i,i+1) = 1$ for $i = 0,1,2,\ldots$.  The graph $\Gamma$ is a weighted chain, with vertices $\{0,1,2,\ldots\}$ and multiplicities $m(i,i+1) = i+1$.  The ``up'' operators in these graphs correspond to multiplication by $x$ in $\Z[x]$ and $\Z[x,x^2/2,x^3/3!,\ldots]$.  This pair of dual graded graphs occurred as Example 2.2.1 in \cite{F}.

\subsection{0-Hecke algebras}
The {\it 0-Hecke algebra} $H_n(0)$ is the unital algebra over $\C$ generated by $T_1, T_2, \ldots, T_{n-1}$ with relations
\begin{align*}
T_i^2 &= -T_i \\
T_iT_j &= T_jT_i  & \mbox{for $|i-j| > 1$} \\
T_i T_{i+1}T_i &= T_{i+1}T_iT_{i+1}.
\end{align*}
The $0$-Hecke algebra has a basis $\{T_w \mid w \in \s_n\}$ and an external multiplication $H_i(0) \otimes H_j(0) \to H_{i+j}(0)$ defined in a manner similar to the nilCoxeter algebras.  An explicit realization of $H_n(0)$ is obtained by the {\it isobaric divided difference operators}.  The representation theory of the tower $H(0) = \bigoplus_{n \geq 0} H_n$(0) was worked out by Krob and Thibon \cite{KT}.

A composition $I = (i_1,i_2,\ldots,i_r)$ of $n$ is a finite sequence of positive integers summing to $n$.
 The algebra $H_n(0)$ is not semi-simple and it has $2^{n-1}$ non-isomorphic simple modules $S_I$ all of dimension 1, as $I$ ranges over the compositions of $n$.  The projective cover $P_I$ of $S_I$ has dimension $\dim(P_I) = \{w \in \s_n \mid \Des(w) = \Des(I):= \{i_1,i_1+i_2,\ldots,i_1 + \cdots + i_{r-1}\}$. It is known that $G(H(0)) = \QSym$, the Hopf algebra of quasi-symmetric functions, and $K(H(0)) = \NSym$, the Hopf algebra of noncommutative symmetric functions.  The class of $[S_I]$ in $\QSym$ is given by the fundamental quasi-symmetric function $F_I \in \QSym$.  The class of $[P_I]$ in $\NSym$ is given by the ribbon Schur function $R_I \in \NSym$.

The graph $\Gamma'$ is an infinite binary tree with vertices of height $n$ identified with compositions of $n$.  There are edges (with multiplicity 1) joining the composition $(i_1,i_2,\ldots, i_r)$ with the compositions $(1, i_1, i_2, \ldots, i_r)$ and with $(i_1 + 1, i_2, \ldots, i_r)$.  The graph $\Gamma$ has edges (multiplicity 1) joining $(i_1,i_2,\ldots, i_r)$ with
\begin{equation*}
\begin{split}
 \{ (i_1,\ldots,i_{j-1}, i_{j} + 1, i_{j+1},\ldots,i_r),(i_1,\ldots,i_{j-1}, k+1, i_{j} - k, i_{j+1},\ldots,i_r)\}
\end{split}
\end{equation*}
for each $j = 1, 2, \ldots, r$ and $k=0,\ldots,i_j-1$.  We reproduce these graphs in Figure \ref{fig:bin} of Section~\ref{ssec:f0hecke} (the edge labels should be ignored for now).  This pair of dual graded graphs occurred as Example 2.3.6 in \cite{F}.

\subsection{Hecke algebras at roots of unity}\label{ssec:hecke}
Let $v \in \C$.  The Hecke algebra $H_n(v)$ is generated by $T_1,T_2,\ldots,T_{n-1}$ with  relations
\begin{align*}
T_i^2 &= (v-1)T_i + v \\
T_iT_j &= T_jT_i  & \mbox{for $|i-j| > 1$} \\
T_i T_{i+1}T_i &= T_{i+1}T_iT_{i+1}.
\end{align*}

The Hecke algebra has a basis $\{T_w \mid w \in \s_n\}$ and an external multiplication $H_i(v) \otimes H_j(v) \to H_{i+j}(v)$ defined in a manner similar to the nilCoxeter algebras.
If $v = 0$, then we recover the $0$-Hecke algebra $H_n(0)$.  If $v = 1$, then we recover the symmetric groups algebras.  If $v$ is neither $0$ nor a root of unity, then the tower $H(v) = \bigoplus_{n \geq 0} H_n(v)$ has representation theory identical to that of $\bigoplus_{n \geq 0} \C\s_n$.  In this case we say that $v$ is generic.

We now let $v = \zeta$ be a primitive $r$-th root of unity, and let
$H(\zeta) = \bigoplus_{n \geq 0} H_n(\zeta)$ denote the corresponding
tower of algebras.  The representation theory of this infinite
family of towers of algebras is not completely understood.  We refer
the reader to \cite{LLT} for the following discussion.  Let
$\J^{(r)} \subset \Sym$ be the ideal in $\Sym$ generated by the
power symmetric functions $p_r, p_{2r}, p_{3r}, \ldots$.  Then the graded dual Hopf algebras $G(H(\zeta)) = \Sym/\J^{(r)}$ and
$K(H(\zeta)) = (\J^{(r)})^\perp$, where $(\J^{(r)})^\perp \subset
\Sym$ is the set of elements annihilated by $\J^{(r)}$ under the
usual pairing of $\Sym$ with itself.  Ariki \cite{Ari}, proving a
conjecture from \cite{LLT}, showed that the symmetric functions
representing the classes of the simple modules, or the projective
indecomposable modules, can be expressed in terms of Schur functions
via the (lower and upper) global bases at $q=1$ of the Fock space
representation of $U_q(\hat{sl}_r)$.

The graded graphs $\Gamma$ and $\Gamma'$ are not known explicitly to
our knowledge, though they have been the subject of much recent
work; see for example \cite{Bru, LLT}.  In particular, these
branching graphs are closely related to the crystal graphs of
quantum affine algebras of type $A$.  It follows from Theorem
\ref{thm:HNLS} that

\begin{cor} \label{C:rootof1}
The branching graph $\Gamma$ for the simple modules, and the branching graph $\Gamma'$ for the projective indecomposable modules of $H(\zeta)$ form a pair of dual graded graphs with differential coefficient $r = 1$.
\end{cor}


\begin{remark}\label{rem:super}
The case $r\ge 2$ is abundant. In particular, for $r=2$ one can
consider towers of super-algebras and super modules. This is how
Seergev~\cite{ANS} constructed the combinatorial Hopf algebra of
$Q$-Schur functions from the tower of Seergev algebras. This is also
how Bergeron, Hivert, and Thibon~\cite{BHT} constructed the
combinatorial Hopf algebra of $\Theta$-peak quasisymmetric functions
from the tower of Hecke-Clifford algebras. Theorem~\ref{dimthm} also
holds for towers of super-algebras; the proof is a direct adaptation
of the one presented here.
\end{remark}

\section{Two parameter Hecke algebras and conjectural classification}\label{sec:HC}
Let $a, b \in \C$.  Let $H_n(a,b)$ denote the {\it two-parameter
Hecke algebra} with generators $T_1,T_2,\ldots,T_{n-1}$ and
relations
\begin{align*}
T_i^2 &= aT_i + b \\
T_iT_j &= T_jT_i  & \mbox{for $|i-j| > 1$} \\
T_i T_{i+1}T_i &= T_{i+1}T_iT_{i+1}.
\end{align*}

\begin{proposition}\label{prop:twoparam}
The $\C$-algebra $H_n(a,b)$ is isomorphic to one of the following four families of algebras:
\begin{enumerate}
\item[H1]
a Hecke algebra $H_n(v)$ at a generic (see subsection \ref{ssec:hecke}) value of $v$, or
\item[H2]
a Hecke algebra $H_n(\zeta)$ at a root of unity $\zeta$, or
\item[H3]
the 0-Hecke algebra $H_n(0)$ (when $a \neq 0$ but $b = 0$), or
\item[H4]
the nilCoxeter algebra $N_n$ (when $a = b = 0$).
\end{enumerate}
\end{proposition}

\begin{proof}
If $(a,b) = (0,0)$, then $H_n(a,b)=N_n$ the nilCoxeter algebra.  Otherwise, we can find a non-zero $z \in \C$ satisfying
\begin{equation}\label{eq:quad}
az = bz^2 - 1.
\end{equation}
The elements $T'_i = zT_i$ then satisfy
$$
(T'_i)^2 = (q-1)T'_i + q
$$
where $q = bz^2$.  Note that the braid relation for the $T_i$ implies the braid relation for the $T'_i$.  Thus
$H_n(a,b)$ is isomorphic to $H_n(bz^2)$.

If $b = 0$, then $H_n(a,0)$ is isomorphic to the 0-Hecke algebra $H_n(0)$ and we are in Case (H3).
Otherwise we are in Case (H1) or (H2).  Note that if $z$ and $z'$ are the two roots of (\ref{eq:quad}), then $bz^2$ is a $r$-th root of unity if and only if $b(z')^2$ is.  This follows from the fact that $zz' = -1/b$.
\end{proof}

Note that the isomorphism of Proposition \ref{prop:twoparam} is compatible with the external multiplication of the obvious construction of the tower $H(a,b) = \bigoplus_{n \geq 0} H_n(a,b)$.
It thus follows that for any $a,b \in \C$, the tower $H(a,b)$ gives rise to one of the graded dual Hopf algebras, and dual graded graphs, discussed in Section \ref{sec:ex}.

Based on this and Theorem \ref{dimthm}, we conjecture

\begin{conjecture} \label{C:class} \mbox{}
\begin{enumerate}
\item (Weak version) Suppose $A$ is a tower of algebras with $\dim(A_1) = 1$, giving rise to graded dual Hopf algebras $K(A)$ and $G(A)$.  Then the pair $(K(A),G(A))$ is isomorphic, together with their distinguished bases (classes of simples and indecomposable projectives), to one of the examples in Section \ref{sec:ex}.
\item (Strong version) Suppose $A$ is a tower of algebras with $\dim(A_1) = 1$, giving rise to graded dual Hopf algebras.  Then $A$ is isomorphic to one of the towers $H(a,b)$.
\end{enumerate}
\end{conjecture}

Zelevinsky \cite{AVZ} shows that a graded connected self-dual Hopf algebra $H$, with a self-dual basis $\{b_\lambda\}$ such that all product and coproduct structure constants are positive with respect to this basis, must be a tensor product of the Hopf algebra of symmetric functions, together with the tensor product of the Schur function basis.  Thus Conjecture \ref{C:class}(1) holds when $A$ is a tower of semisimple algebras.

\section{Quantum version}\label{Sec:Q}
In this section, we describe a ``quantum'' version of our theorem.  We replace \eqref{eq:comp} with
$$
\text{tower of filtered algebras} \rightarrow \text{$q$-twisted Hopf
algebra} \rightarrow \text{$q$-dual graded graph}
$$
We shall not make the first arrow completely axiomatic here. 

\subsection{From filtered towers of algebras to $q$-twisted Hopf algebras}\label{ssec:filtered}
We first recall the notion of a $q$-twisted Hopf algebra \cite{LZ,HNT2}.  Let $H$ be a graded connected algebra over $\Z[q]$, equipped with an associative graded coproduct $\Delta: H \to H \otimes_{\Z[q]} H$.  The formula for the \textit{$q$-twisted product of tensors} is
\begin{equation*}
(a\otimes b)\cdot_q (a'\otimes b')=q^{{\rm deg}(b)\cdot {\rm deg}(a')}(aa'\otimes bb').
\end{equation*}
We say that $H$ is a {\it $q$-twisted Hopf algebra} if $\Delta(a)\cdot_q\Delta(b) = \Delta(ab)$ for every $a, b \in H$.  The other structure maps (unit, counit, antipode) will not concern us here.

The notion of $q$-twisted Hopf algebras is a particular instance of {\it twisted} Hopf algebras (or Hopf monoids) for braided monoidal categories~\cite{twistHopf}. Indeed, a $q$-twisted Hopf algebra is a Hopf monoid in the category of graded vector spaces with braiding induced by $\tau_q\colon V \otimes W \to W \otimes V$ where $\tau_q(v\otimes w)=q^{nm} w\otimes v$ for $v$ (resp. $w$) a homogeneous element of degree $n$ (resp. $m$).

Now let $A = \bigoplus_{n \geq 0} A_n$ be a tower of algebras as in Section \ref{sec:BL}.  We suppose that each $A_n$ is equipped with a filtration $A_n = A_n^{(0)} \supset A_n^{(1)} \supset \cdots $ such that each $A_n^{(k)}$ is a left ideal in $A_n$.  We call this a {\it filtered tower of algebras}.  Let $M$ be a left $A_n$-module and $M' \subset M$ be a subset.  If $M=A_n\cdot M'$, then the sequence
$M^{(0)} = A_n^{(0)} \cdot M' \supset M^{(1)} = A_n^{(1)} \cdot M' \supset M^{(2)} = A_n^{(2)} \cdot M' \supset \cdots$
is a filtration of $M$ by left submodules of $A_n$.  The {\it graded character} $[M]_q \in G_0(A) \otimes_\Z \Z[q]$ is defined by
$$
[M]_q = \sum_{i \geq 0} q^i \; [M^{(i)}/M^{(i+1)}].
$$
Obviously $[M]_q$ depends on $M'$ even though it is suppressed in the notation.

We now define a multiplication $*$ in $G(A)_q = G(A) \otimes_\Z \Z[q]$.  For
$[M]\in G_0(A_m)$ and $[N]\in G_0(A_n)$ we let
 $$[M]*[N]=\left[\mbox{Ind}^{A_{m+n}}_{A_{m}\otimes
A_{n}}M\otimes N\right]_q $$
with respect to the subset $M \otimes N \subset \mbox{Ind}^{A_{m+n}}_{A_{m}\otimes
A_{n}}M\otimes N$.  We also equip $G(A)_q$ with the usual coproduct of $G(A)$, extended by linearity to $G(A)_q$.  Assume:

\begin{enumerate}
\item[(Q1)] The multiplication $*$ is a well-defined associative product on $G(A)_q$.
\item[(Q2)] $G(A)_q$ is a $q$-twisted Hopf algebra.
\item[(Q3)] Graded characters and inductions are sufficiently compatible such that
$$\beta^{*n} = [A_n]_q = \sum_{i \geq 0} q^i \; [A_n^{(i)}/A_n^{(i+1)}]$$
where $\beta$ is as in Section \ref{sec:proof}.
\end{enumerate}

\begin{remark}
For example, (Q3) would follow from the more general compatibility equation
$$
[M_1]*[M_2]*\cdots*[M_r] = \left[\mbox{Ind}^{A_{m_1+m_2+\cdots+m_r}}_{A_{m_1}\otimes \cdots \otimes
A_{m_r}}M_1 \otimes \cdots \otimes M_r \right]_q
$$
and some assumption about the structure map $A_1 \otimes \cdots \otimes A_1 \to A_n$.  We believe that this compatibility relation is the most natural one, but do not need such generality here.
\end{remark}

For our purposes, the precise construction of $G(A)_q$ is not crucial, as long as (Q1)-(Q3) are satisfied.  We say that a multiplication on $G(A) \otimes_\Z \Z[q]$ quantizes $G(A)$ if it reduces to the usual multiplication of $G(A)$ at $q=1$.  Let $[r] = 1+ q + \cdots +q^{r-1}$, and $[r]! = [r][r-1]\cdots[1]$ be the usual $q$-analogues.

\begin{theorem}\label{T:quantization}
Let $A$ be a filtered tower of algebras, and suppose a quantization $G(A)_q$ of $G(A)$ exists, satisfying (Q1)-(Q3) above.  Then
$$
\dim_q(A_n) = \sum_{i \geq 0} q^i \, \dim(A^{(i)}_n/A^{(i+1)}_n)= r^n [n]!
$$
where $r = \dim(A_1)$.
\end{theorem}
Theorem \ref{T:quantization} will be proved in subsection \ref{ssec:qproof} below.

\subsection{From $q$-twisted Hopf algebras to quantized dual graded graphs}
Quantized dual graded graphs are defined and studied in \cite{Lam}.  We now allow our graded graphs $\Gamma = (V, E, h,m)$ to have multiplicities taking values in $\N[q]$, where $\N = \{0,1,2,\ldots\}$ (for some purposes $\N[q^{1/2},q^{-1/2}]$ could also be considered).  Making definitions analogous to those in Section \ref{sec:HNLS}, we say that $(\Gamma,\Gamma')$ is a pair of {\it quantized dual graded graphs} with differential coefficient $r$, if the linear operators $U_\Gamma$ and $D_{\Gamma'}$ satisfy
$$
D_{\Gamma'}U_\Gamma - q\, U_\Gamma D_{\Gamma'} = r {\rm Id}.
$$
For now we allow $r$ to lie in $\Z[q]$, though in the end we have no need for such generality.

We define $f^v_\Gamma \in \N[q]$ as before.  The following is the quantized analogue of Theorem \ref{thm:F}.

\begin{theorem}[Lam~\cite{Lam}] \label{thm:Lam} Let $(\Gamma,\Gamma')$ be a pair of quantized dual graded graphs
with differential coefficient $r$. Then
$$r^n [n]!=\sum_{v \colon h(v)=n}f^v_{\Gamma}f^v_{\Gamma'}.$$
\end{theorem}

We now generalize Theorem \ref{thm:HNLS} to the quantized setting.  We first make the general observation that the graded dual of a graded $q$-twisted Hopf algebra over $\Z[q]$ is again a graded $q$-twisted Hopf algebra.

Let $H_\bullet = \bigoplus_{n \geq 0} H_n$ and $H^\bullet =
\bigoplus_{n \geq 0} H^n$ be graded dual $q$-twisted Hopf algebras over $\Z[q]$
with respect to the pairing $\scal{\,.\,,.\,}: H_\bullet \times H^\bullet
\to \Z[q]$. We assume that we are given dual sets of homogeneous free
$\Z[q]$-module generators $\{p_\lambda \in H_\bullet\}_{\lambda \in
\Lambda}$ and $\{s_\lambda \in H^\bullet\}_{\lambda \in \Lambda}$,
such that all structure constants lie in $\N[q]$.  We
also assume that $\dim(H_i) = \dim(H^i) < \infty$ for each $i \geq
0$ and $\dim(H_0) = \dim(H^0) = 1$, so that $H_0$ and $H^0$ are
spanned by distinguished elements the unit $1$.  Let us suppose we
are given non-zero homogeneous elements $\alpha \in H_1$ and $\beta
\in H^1$ of degree 1 such that $\alpha p_{\mu}$ (resp. $\beta s_{\mu}$)
is a linear combination of $\{p_{\lambda}\}$ (resp. $\{s_{\lambda}\}$) with $\N[q]$-coefficients for any $\mu\in\Lambda$.

We now define graded graphs $\Gamma(\beta)$ and $\Gamma'(\alpha)$ exactly as in Section \ref{sec:HNLS}.  The following theorem generalizes Theorem \ref{thm:HNLS}.
%

\begin{theorem}\label{thm:qHNLS}
The graded graphs $\Gamma = \Gamma(\beta)$ and $\Gamma' =
\Gamma'(\alpha)$ form a pair of quantized dual graded graphs with {\it
differential coefficient} $\scal{\alpha,\beta}$.
\end{theorem}
\begin{proof}
The proof is identical to that of Theorem \ref{thm:HNLS} until the final calculation, which proceeds
\begin{align*}
D_{\Gamma'} U_{\Gamma}(x) & = D_{\Gamma'}(\beta x)\\
&=\sum_{\mu \in \Lambda} \scal{\alpha \otimes p_\mu, \Delta(\beta)._q\Delta(x)}s_\mu \\
&=\sum \scal{\alpha \otimes p_\mu, \beta \, x^{(1)} \otimes x^{(2)}+ q^{{\rm deg}(x^{(1)})} x^{(1)} \otimes \beta \, x^{(2)}} s_\mu\\
& = \sum \left(\scal{\alpha,\beta \, x^{(1)}}\,x^{(2)} + q^{{\rm deg}(x^{(1)})} \, \scal{\alpha,x^{(1)}}\,\beta\,x^{(2)}\right) \\
&= \scal{\alpha,\beta} x + \sum_{{\rm deg}(x^{(1)}) = 1} q \, \beta \, \scal{\alpha,x^{(1)}}\, x^{(2)}\\
&= \scal{\alpha,\beta}x +qU_{\Gamma}D_{\Gamma'}(x).
\end{align*}
\end{proof}

\subsection{Proof of Theorem \ref{T:quantization}}
\label{ssec:qproof}
The proof is analogous to that of Theorem \ref{dimthm}.  By Theorems \ref{thm:Lam} and \ref{thm:qHNLS}, and assumptions (Q1)-(Q2), it suffices to show that $\dim_q(A_n) = \sum_{i}f^{s_i^{(n)}}_{\Gamma}f^{s_i^{(n)}}_{\Gamma'}$.  By assumption (Q3),
$$
 [A_n]_q = \beta^{*n} = \sum f^{s_i^{(n)}}_{\Gamma} [S_i^{(n)}]
$$
in $G(A)_q$.  But comultiplication in $G(A)_q$ is the same as in $G(A)$, so by Lemma \ref{lem:dimpath}, we have
$$
\dim_q(A_n) = \sum_i f^{s_i^{(n)}}_{\Gamma} \dim(S_i^{(n)}) = \sum_{i}f^{s_i^{(n)}}_{\Gamma}f^{s_i^{(n)}}_{\Gamma'},
$$
as required.

\subsection{Filtered nilCoxeter algebras}
Recall the nilCoxeter algebra $N_n$ from subsection \ref{ssec:nilCox}.  The algebras $N_n$ are filtered by $N_n^{(k)}=\oplus_{\ell(w)\geq k}\mathbb{C}T_w$, where $\ell(w)$ denotes the length of the permutation $w$.  In fact, $N_n$ is a graded algebra, and the graded representation theory was also considered by Khovanov \cite{Kho}.  (The formulae for our $q$-twisted Hopf algebra below differs somewhat from that considered in \cite{Kho}.)

The construction of subsection \ref{ssec:filtered} produces a
$q$-twisted Hopf algebra $G(N)_q$ with multiplication as follow.
Let $S_n$ and $S_m$ be the (unique) simple modules of $N_n$ and $N_m$
respectively. Then
    $$\left[\mbox{Ind}^{N_{m+n}}_{N_{m}\otimes N_{n}}S_m\otimes S_n\right]_q = {[n+m]! \over [n]![m]!} [S_{n+m}].$$
We can thus identify $G(N)_q$  with $\Z[q][x/[1], x^2/[2]!, x^3/[3]!, \cdots]$, equipped with the usual multiplicative structure. The coproduct is then defined by
$$
\Delta\left(\frac{x^c}{[c]!}\right) = \sum_{r= 0}^c \frac{x^r}{[r]!} \otimes \frac{x^{c-r}}{[c-r]!}.
$$
The $q$-twisted structure reduces to a well-known identity for $q$-binomial coefficients.  Let $\binom{m}{n}_q = [m]!/([n]![m-n]!)$ be the usual $q$-binomial coefficients, where by convention $\binom{m}{n}_q = 0$ if $m < n$.  Then the following identity is standard (and follows easily from the interpretation of $\binom{m}{n}_q$ as the rank-generating function of the product of two chains):
\begin{equation}\label{E:binom}
\binom{m}{n}_q = \sum_{i=0}^n q^{(n-i)(r-i)}\,\binom{r}{i}_q\,\binom{m-r}{n-i}_q
\end{equation}
for any $1 \leq r \leq n$.  We then calculate
 \begin{align*}
&\Delta\left(\frac{x^a}{[a]!}\right)._q \Delta\left(\frac{x^{c-a}}{[c-a]!}\right)\\
&=\left(\sum_{i= 0}^a \frac{x^i}{[i]!} \otimes \frac{x^{a-i}}{[a-i]!}\right)._q\left(\sum_{j= 0}^{c-a}\frac{x^j}{[j]!} \otimes \frac{x^{c-a-j}}{[c-a-j]!}\right)\\
&= \sum_{i=0}^a \sum_{j=0}^{c-a} q^{(a-i)j}\, \binom{i+j}{i}_q \, \binom{c-i-j}{a-i}_q \, \frac{x^{i+j}}{[i+j]!} \otimes \frac{x^{c-i-j}}{[c-i-j]!}.
\end{align*}
Now take the coefficient of $x^r/[r]! \otimes x^{c-r}/[c-r]!$ and use \eqref{E:binom} with $m = c$ and $n = a$ to see that this is equal to $\binom{c}{a}_q \, \Delta(x^c/[c]!)$.  Theorem \ref{T:quantization} reduces to the well-known combinatorial identity $\sum_{w \in \s_n} q^{\ell(w)} = [n]!$.

The graph $\Gamma'$ is still a chain, as in subsection \ref{ssec:nilCox}.  The graph $\Gamma$ has edge multiplicities $m(i,i+1) = [i+1]$.  The pair $(\Gamma,\Gamma')$ is a pair of quantized dual graded graphs with differential coefficient $r = 1$.


\subsection{Filtered 0-Hecke algebras}\label{ssec:f0hecke}
\def\Shuf{{\rm Shuf}}

\begin{figure}
\begin{center}
\begin{picture}(20,150)(150,-10)
\put(54,-10){$\Gamma'(\alpha)$}
\put(60,5){$\emptyset$}
\put(60,35){1}
\put(20,65){2}\put(100,65){11}

\put(0,95){3}
\put(40,95){12}\put(80,95){21}\put(120,95){111}

\put(-10,135){4}\put(10,135){ 13}
\put(30,135){22}\put(50,135){112}\put(70,135){31}\put(90,135){121}\put(110,135){211}\put(130,135){1111}

\put(62,14){\line(0,1){17}}
\put(67,42){\line(3,2){30}}
\put(56,42){\line(-3,2){30}}
\put(110,72){\line(1,2){10}}
\put(96,72){\line(-1,2){10}}
\put(30,72){\line(1,2){10}}
\put(16,72){\line(-1,2){10}}
\put(126,104){\line(1,4){7}}
\put(122,104){\line(-1,4){7}}
\put(86,104){\line(1,4){7}}
\put(82,104){\line(-1,4){7}}
\put(45,104){\line(1,4){7}}
\put(41,104){\line(-1,4){7}}
\put(4,104){\line(1,4){7}}
\put(0,104){\line(-1,4){7}}

\put(234,-10){$\Gamma(\beta)$}
\put(241,5){$\emptyset$}
\put(240,35){1}
\put(200,65){2}\put(280,65){11}

\put(180,95){3}
\put(220,95){12}\put(260,95){21}\put(300,95){111}

\put(170,135){4}\put(186,135){13}
\put(210,135){22}\put(230,135){112}\put(250,135){31}\put(270,135){121}\put(290,135){211}\put(310,135){1111}

\put(242,14){\line(0,1){17}}
\put(247,42){\line(3,2){30}}
\put(236,42){\line(-3,2){30}}
\put(290,72){\line(1,2){10}}
\put(276,72){\line(-1,2){10}}
\put(210,72){\line(1,2){10}}
\put(196,72){\line(-1,2){10}}
\put(306,104){\line(1,4){7}}
\put(302,104){\line(-1,4){7}}
\put(265,104){\line(1,4){7}}
\put(261,104){\line(-1,4){7}}
\put(225,104){\line(1,4){7}}
\put(221,104){\line(-1,4){7}}
\put(184,104){\line(1,4){7}}
\put(180,104){\line(-1,4){7}}

\put(214,72){\line(5,2){48}}
\put(272,72){\line(-5,2){48}}

\put(186,104){\line(1,1){27}}
\put(189,104){\line(2,1){59}}

\put(300,104){\line(-1,1){27}}
\put(297,104){\line(-2,1){59}}

\put(266,104){\line(1,1){27}}
\put(260,104){\line(-3,2){43}}

\put(220,104){\line(-1,1){27}}
\put(226,104){\line(3,2){43}}

{\tiny\color{blue}\put(237,19){1}
\put(264,50){$q$}
\put(217,48){1}

\put(294,76){$q^2$}
\put(275,75){1}
\put(264,77){$q$}

\put(207,76){$q$}
\put(196,74){1}
\put(218,77){$q^2$}

\put(308,109){$q^3$}
\put(300,111){1}
\put(294,111){$q$}
\put(284,105){$q^2$}

\put(262,110){$q$}
\put(255,109){1}
\put(270,105){$q^3$}
\put(244,105){$q^2$}

\put(227,109){$q^2$}
\put(220,109){1}
\put(211,106){$q$}
\put(232,104){$q^3$}

\put(180,109){$q$}
\put(174,107){1}
\put(187,109){$q^2$}
\put(196,104){$q^3$}
}

\end{picture}
\end{center}
\caption{Quantizations of the dual graded graphs for the $0$-Hecke algebra}
\label{fig:bin}
\end{figure}
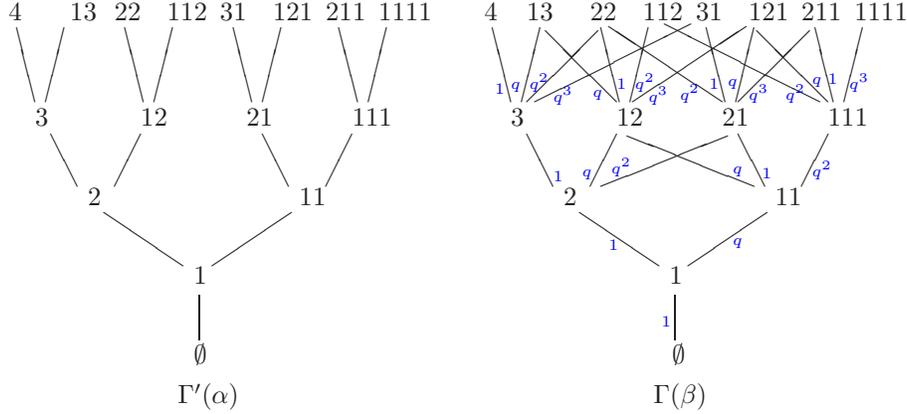

Now consider the tower $H(0)$ of $0$-Hecke algebras, which we equip
with the filtrations $H_n(0)^{(k)}=\oplus_{\ell(w)\geq
k}\mathbb{C}T_w$, where $\ell(w)$ denotes the length of the
permutation $w$.  The representation theory of this filtered tower
of algebras was considered by Thibon and Ung \cite{TU}, to which we
refer for further details.

The construction of subsection \ref{ssec:filtered} produces a $q$-twisted Hopf algebra $G(H(0))_q = \QSym_q$ known as the {\it quantum quasi-symmetric functions}.  $\QSym_q$ is spanned over $\Z[q]$ by the fundamental quasi-symmetric functions $F_I$, equipped with the {\it quantum shuffle product}, which we now describe.  Let $I = (i_1,i_2,\ldots, i_r)$ and $J = (j_1,j_2,\ldots,j_s)$ be compositions of $m$ and $n$ respectively.  Let $w = w_1w_2 \cdots w_m \in \s_m$ and $v = v_1 v_2 \cdots v_n \in \s_n$ be permutations with descent sets $\Des(w) = \Des(I):=\{i_1, i_1 + i_2, \ldots, i_1 + i_2 + \cdots + i_{r-1}\}$ and $\Des(v) = \Des(J)$.   Denote the shuffles of $w$ and $v$ by $\Shuf(w,v) \subset \s_{m+n}$, which we illustrate with an example: if $w = 132$ and $v = 21$, then $$\Shuf(w,v) = \{13254, 13524,15324,51324,13542,15342,51342,15432,51432,54132\}.$$  For $u \in \s_n$, write $C(u)$ for the composition $C$ of $n$ such that $\Des(C) = \Des(u)$.  Then in $\QSym_q$,
$$
F_I * F_J = \sum_{u \in \Shuf(w,v)} q^{\theta(u)}\, F_{C(u)}
$$
where
$$
\theta(u) = \#\{(i,j) \in \{1,2,\ldots,m\} \times \{m+1,\ldots,m+n\} \mid \mbox{$i$ occurs after $j$ in $u$}\}.
$$
With this product, and the usual coproduct of quasi-symmetric functions, $\QSym_q$ becomes a $q$-twisted Hopf algebra.  At $q = 1$, the algebra $\QSym_q$ reduces to $\QSym$.  We caution that while $\QSym$ is commutative, the algebra $\QSym_q$ is noncommutative, and in fact is isomorphic to $\NSym$ as a ring (\cite{TU}).  Theorem \ref{T:quantization} reduces, again, to the well-known combinatorial identity $\sum_{w \in \s_n} q^{\ell(w)} = [n]!$.

The graph $\Gamma'$ is the infinite lifted binary tree as before.  The edge multiplicities for $\Gamma$ are powers of $q$, which is illustrated in Figure \ref{fig:bin}.  
The edge joining $(i_1,i_2,\ldots, i_r)$ with $(i_1,\ldots,i_{j-1}, i_{j} + 1, i_{j+1},\ldots,i_r)$ has multiplicity $q^{i_1 + i_2 + \cdots + i_{j-1}}$; the edge joining $(i_1,i_2,\ldots, i_r)$ with $ (i_1,\ldots,i_{j-1},k+1, i_{j} - k, i_{j+1},\ldots,i_r)$ has multiplicity $q^{i_1 + \cdots + i_{j-1} + k + 1}$.  This gives a pair of quantized dual graded graphs with differential coefficient $r = 1$, which is not difficult to verify directly.

\section{Further directions}
\label{sec:further} There are many new avenues that could be
explored from the point of view developed in this paper.

\subsection{Generalized bialgebras, and Hopf monoids}
In this paper, we have three worlds connected by some constructions:
\begin{equation}\label{eq:compagain}
\text{tower of algebras} \longrightarrow \text{combinatorial Hopf
algebra} \longrightarrow \text{dual graded graph.}
\end{equation}
It is natural to ask if this is possible for other triples of
similar objects.  In particular, in \cite{BL}, Bergeron and Li
described how a general tower of algebras gives rise to generalized
bialgebras in the sense of Loday~\cite{Loday}.  Fomin's notion of
dual graded graphs is naturally related to Hopf algebras.  It is
thus natural to ask what generalization of dual graded graph is
related to other generalized bialgebras.  In particular, what
relation replaces $DU- UD = rI$?

The theory developed in Section~\ref{Sec:Q} suggests another
direction.  The construction there replaces Hopf algebras with Hopf
monoids which lie in a different braided monoidal category. It would
be interesting to study if it is possible to construct a triple for
other types of Hopf monoids.

\subsection{Category of Combinatorial Hopf Algebras}

In~\cite{ABS}, Aguiar, Bergeron, and Sottile considered the category
of combinatorial Hopf algebras, consisting of pairs $(H,\zeta)$
where $H$ is a graded connected Hopf algebra and $\zeta\colon H\to
\C$ is a character (multiplicative homomorphism to the ground field
$\C$).

The first arrow of (\ref{eq:compagain}) allows us to construct very
natural pairs $(H,\zeta)$. More precisely, let $A=\bigoplus_{n\ge 0}
A_n$ be a tower of algebras satisfying Theorem~\ref{dimthm} and
suppose we are given a family $\{P_n^0\}$ of one-dimensional
projective modules, satisfying the following compatibility relation:
\begin{equation}\label{eq:grouplike}
\hbox{Res}^{A_{n+m}}_{A_n\otimes A_m} P_{m+n}^0 = P_n^0\otimes
P_m^0.\end{equation}


We can thus define the linear map
 $$ \begin{array} {rrcl}
  \zeta^0\colon &  G_0(A) & \to & \Z\cr
   &[M]\in G_0(A_n) &\mapsto & \langle [P_n^0], [M]\rangle.\cr
  \end{array}$$
Using (\ref{eq:grouplike}), it is clear that $\zeta^0$ is a
character (taking values on $\Z$). We thus have that
$(G_0(A),\zeta^0)$ is a combinatorial Hopf algebra in the sense of
\cite{ABS}.

In our key examples, $(G_0(A),\zeta^0)$ satisfies some universal
properties. For instance, for the tower of $0$-Hecke algebras, a
natural family of one-dimensional projective modules exists and in
$\NSym$ are encoded by the ribbon Schur functions $R_{(n)}$. The
resulting character $\zeta^0$ is precisely what is needed to get the
universal property of $\QSym$ as in~\cite{ABS}. For the tower of
$\s_n$-group algebras, the one-dimensional projective modules are
encoded by the Schur functions $s_{(n)}$ in $\Sym$. As shown
in~\cite{ABS}, $(\Sym,\zeta^0)$ is universal among cocomutative Hopf
algebras.

In both cases, this family of one-dimensional projective modules are
``trivial'' -- both symmetric group algebras and $0$-Hecke algebras
have a distinguished basis (see Section \ref{sec:ex}) which acts
trivially on these modules.  It is tempting to say that as soon as a
tower has ``trivial'' projective modules, then it is a universal
object in some category. This seems to be the case for most of the
examples we know, and P. Choquette \footnote{Personal communication with P. Choquette, 2009.} has some results along
this line for our quantum example in Section~\ref{Sec:Q}. Yet we do
not have such a result for the tower of Hecke algebras at root of
unity. It would be very interesting to find a category for which
$\Sym/\mathcal{J}^{(r)}$ and its character is a universal object.

\begin{remark}
The towers given in Remark~\ref{rem:super} also have a natural
compatible family of one-dimensional projective modules.  The Hopf
algebra of $Q$-Schur functions (with the associated character) is
the universal object in the category of cocomutative odd combinatorial
Hopf algebras. The Hopf algebra of peak quasisymetric functions is
universal in the category of  odd combinatorial Hopf algebras.
\end{remark}

 \begin{remark}
The tower of nilCoxeter algebras does not have one-dimensional
projective modules as described above.  Yet, the simple modules
$S_n^0$ do satisfy (\ref{eq:grouplike}), giving rise to a character
$\zeta^0$ on $K(A) =\Z[x]$, given by $\zeta^0(x^n) = 1$. The
universal property satisfied by $\mathbb{Z}[x]$ is somehow trivial:
for any combinatorial Hopf algebra $(H,\zeta)$, with $\zeta$ taking
values in $\Z$, there is a unique algebra morphism
$\hat{\zeta}\colon H\to\Z[x]$ defined by
$\hat{\zeta}(h)=\zeta(h)x^n$ for $h\in H_n$, satisfying $\zeta^0
\circ \hat\zeta = \zeta$.
\end{remark}

\subsection{Bi-tower of algebras, and categorification}
One can consider Fomin's work on dual graded graphs \cite{F} as
generalizing Stanley's notion of differential posets by considering
different posets for the up and for the down operators.  In this
context it seems natural to allow two distinct tower structures on a
family of algebras in order to define induction and restriction with
a compatibility relation.  More precisely, let us say that a
bi-tower of algebras $(A,\rho,\rho')$ is a tower
$A=\bigoplus_{n\ge0} A_n$ such that $\rho\colon A\otimes A\to A$ is
a tower of algebras satisfying Conditions (1), (2) and (3) of
Section~\ref{sec:BL}, and $A$ with $\rho'\colon A\otimes A\to A$
also satisfies (1), (2) and (3).  Now, we use $\rho$ to define the
multiplication of $G(A)$ and the comultiplication of $K(A)$ but we
use $\rho'$ to define the comultiplication of $G(A)$ and the
multiplication of $K(A)$.

We may now ask whether $(K(A),G(A))$ form a pair of dual graded Hopf
algebras.  It is straightforward to check that Theorem~\ref{dimthm}
also holds for bi-towers of algebras which give rise to a pair of
graded dual Hopf algebras.

\begin{theorem}
Let $(A,\rho,\rho')$ be a bi-tower of algebras such that its
associated Grothendieck groups form a pair of graded dual Hopf
algebras.  Then $\dim(A_n)=r^nn!$ where $r=\dim(A_1)$.
\end{theorem}

\begin{example}\label{ex:fund}
We give one interesting example of a bi-tower of algebras. Let
$A_n=\mathbb{C}^{n!}$ be the commutative semisimple algebra of
dimension $n!$. This implies that the Grothendieck groups
$G(A)=K(A)=A$. The canonical basis of $A_n$ is given by
$\{e_\sigma:\sigma\in \s_n\}$. To define $\rho$ and $\rho'$, let
$\hbox{Shuf}(n,m)=\{\zeta\in \s_{n+m}: \zeta(1)<\cdots <\zeta(n),
\zeta(n+1)<\cdots<\zeta(n+m)\}$. We also consider the canonical
imbedding $\s_n\times \s_m\hookrightarrow \s_{n+m}$ and denote by
$\sigma\times\pi\in \s_{n+m}$ the element corresponding to
$(\sigma,\pi)\in \s_n\times \s_m$. We define $\rho\colon A\otimes A\to
A$ by
  $$\rho(e_\sigma\otimes e_\pi)=\sum_{\zeta\in \hbox{Shuf}(n,m)} e_{(\sigma\times\pi)\zeta^{-1}}.$$
  We define $\rho'\colon A\otimes A\to A$ by
  $$\rho'(e_\sigma\otimes e_\pi)=\sum_{\zeta\in \hbox{Shuf}(n,m)} e_{\zeta(\sigma\times\pi)}.$$
 With these two maps, the reader can easily verify that $G(A)$ and its dual $K(A)$ will correspond to the two descriptions of the Malvenuto-Reutenauer Hopf algebras as in~\cite{AS}. The pair of dual graphs corresponding to it is the fundamental pair given in Section 2.6 of~\cite{F}.
\end{example}

We now want to see the diagram in (\ref{eq:compagain}) as functors
between categories.
\begin{equation}\label{eq:compCAT}
{\mathcal T} \longrightarrow {\mathcal H} \longrightarrow {\mathcal G}
\end{equation}
It is very natural to allow bi-towers of algebras in our
constructions. The objects in the first category are bi-towers of
algebras $(A,\rho,\rho')$ which give rise to graded dual Hopf
algebras.  A morphism $F\colon (A,\rho,\rho')\to
(B,\overline{\rho},\overline{\rho}')$ is given by a family of
algebra homomorphisms $F_n\colon A_n\to B_n$ such that $\overline{\rho}\circ
(F\otimes F)=F\circ\rho$ and $\overline{\rho}'\circ (F\otimes
F)=F\circ\rho'$. Moreover we require that for every
primitive idempotent $g$ of $A_{n+m}$, we can find idempotents $e$'s
and $f$'s such that
  $$ gA_{n+m} \cong \bigoplus eA_n\otimes fA_m \qquad\hbox{ and } \qquad F(g)B_{n+m} \cong \bigoplus F(e)B_n\otimes F(f)B_m .$$
As proven by Li in~\cite {Li}, this induces graded dual Hopf algebra
morphisms $F_*\colon K(A)\to K(B)$ and $F^*\colon G(B)\to G(A)$.

We now consider the category $\mathcal H$ with objects
$(H_\bullet,\{p_\lambda\},\alpha,\beta)$ where $H_\bullet=\bigoplus
H_n$ is a graded connected Hopf algebra over $\mathbb Z$, the set
$\{p_\lambda\}$ is a homogeneous basis of $H_\bullet$ such that all
structure constants for multiplication and comultiplication are
non-negative and $\alpha\in H_1$ is a non-negative $\mathbb
Z$-linear combination of the basis  $\{p_\lambda\}$. We denote by
$H^\bullet$ the graded dual of $H_\bullet$ and by  $\{s_\lambda\}$
the homogeneous basis dual to  $\{p_\lambda\}$. The element
$\beta\in H^1$ is a non-negative $\mathbb Z$-linear combination of
the basis  $\{s_\lambda\}$. A morphism $T\colon
(H_\bullet,\{p_\lambda\},\alpha,\beta)\to
(\overline{H}_\bullet,\{\overline{p}_\lambda\},\overline{\alpha},\overline{\beta})$
in this category corresponds to a graded Hopf algebra morphism
$T_\bullet \colon H_\bullet\to \overline{H}_\bullet$ such that
$T_\bullet(p_\lambda)$ is a non-negative linear combination of the
$\{\overline{p}_\lambda\}$ and
$T_\bullet(\alpha)=\overline{\alpha}$. By duality this induces a
graded Hopf algebra morphism $T^\bullet\colon
\overline{H}^\bullet\to H^\bullet$ for which
$T^\bullet(\overline{s}_\lambda)$  is a non-negative linear
combination of the $\{s_\lambda\}$. We require that
$T^\bullet(\overline{\beta})=\beta$.

From the above construction and Section~\ref{sec:BL}, if $F\colon
(A,\rho,\rho')\to (\overline{A},\overline{\rho},\overline{\rho}')$
is a morphism of bi-towers of algebras, then $F_*\colon K(A)\to
K(\overline{A})$ is a graded Hopf algebra morphism such that
$F_*(p_\lambda)$ decomposes into a non-negative linear combination
of projective $\overline{A}$-modules, which is a non-negative linear
combination of the $\{\overline{p}_\lambda\}$. Using the $\alpha\in
K(A_1)$ and $\beta\in G(A_1)$ as defined in Section~\ref{sec:proof},
we obtain that $F_*(\alpha)=\overline{\alpha}$ and
$F^*(\overline{\beta})=\beta$. Hence the construction
$(A,\rho,\rho')\mapsto(K(A),\{p_\lambda\},\alpha,\beta)$ is
functorial.

The third category $\mathcal G$ consists of dual graded graphs $(\Gamma,\Gamma')$. A morphism $\varphi\colon (\overline\Gamma,\overline\Gamma') \to(\Gamma,\Gamma')$ is a $\mathbb Z$-linear map $\varphi\colon {\mathbb Z}\overline{V}\to {\mathbb Z} V$ on the $ {\mathbb Z}$-module of vertices such that $h\circ\varphi=\overline{h}$ where $h$ and $\overline{h}$ are extended linearly, $U_{\Gamma}\circ\varphi =\varphi\circ U_{\overline\Gamma}$ and $D_{\Gamma'}\circ\varphi =\varphi\circ D_{\overline\Gamma'}$. We also require that $\varphi(\overline{v})$ is a non-negative linear combination of ${V}$ for all $\overline{v}\in \overline{V}$.

Given a morphism  $T\colon (H_\bullet,\{p_\lambda\},\alpha,\beta)\to (\overline{H}_\bullet,\{\overline{p}_\lambda\},\overline{\alpha},\overline{\beta})$ in the category $\mathcal H$, we obtain a morphism of the category $\mathcal G$ as follows. First we remark that if $T_\bullet(p_\lambda)=\sum_\mu c_{\lambda,\mu}\overline{p}_\mu$,
then
$T^\bullet\colon \overline{H}^\bullet \to {H}^\bullet$ is a graded $\mathbb Z$-linear map such that $T^\bullet(\overline{s}_\mu)=\sum_\lambda c_{\lambda,\mu}s_\lambda$  is a non-negative linear combination. From Section~\ref{sec:HNLS},
$$U_{\Gamma(\beta)}\circ T^\bullet(\overline{x})=\beta T^\bullet(\overline{x}) = T^\bullet(\overline{\beta}\overline{x})
     =T^\bullet \circ U_{\Gamma(\overline{\beta})}(\overline{x}).
$$
Also
  \begin{align*}
     D_{\Gamma(\alpha)}\circ T^\bullet(\overline{x}) &=  \sum_\lambda \langle \alpha p_\lambda, T^\bullet(\overline{x})\rangle s_\lambda
               =  \sum_\lambda \langle T_\bullet(\alpha p_\lambda), \overline{x}\rangle s_\lambda \cr
               &=  \sum_{\lambda} \sum_\mu c_{\lambda,\mu} \langle \overline{\alpha} \,\overline{p}_\mu, \overline{x}\rangle s_\lambda
               = T^\bullet\Big( \sum_\mu  \langle \overline{\alpha} \,\overline{p}_\mu, \overline{x}\rangle \overline{s}_\mu\Big) \cr
       &=T^\bullet \circ D_{\Gamma(\overline{\alpha})}(\overline{x}).
  \end{align*}
Hence the construction
$(H_\bullet,\{p_\lambda\},\alpha,\beta)\mapsto
(\Gamma(\beta),\Gamma(\alpha))$ is a (contravariant) functor from
$\mathcal H$ to $\mathcal G$. We have thus shown the following
theorem.

\begin{theorem}\label{thm:categorifying}
The two constructions ${\mathcal T}\to{\mathcal H}$ and ${\mathcal H}\to{\mathcal G}$ are functorial.
\end{theorem}

 \begin{remark}
For $r=1$, it should be possible to map the minimal idempotents of a
bi-tower $(A,\rho,\rho')$ (giving rise to graded dual Hopf algebras)
into the bi-tower of Example~\ref{ex:fund} in such a way that we get
a morphism. This would be a good way to see the fundamental role
played by the dual graded graphs given in Section 2.6 of~\cite{F}.
It also would explain the importance of the Malvenuto-Reutenauer
Hopf algebra. This is conceptually plausible but in practice likely
to be very hard. For example, such a morphism is not known for the
tower of $0$-Hecke algebras.
\end{remark}

\end{document}